\documentclass[11pt,a4paper]{article}
\usepackage[latin1]{inputenc}
\usepackage{amsmath}
\usepackage{graphicx}

\title{Solutions to the generalized Towers of Hanoi problem}
\author{Mikael Erik Jörgensen}

\begin{document}


\begin{center}
{\Large Solutions to the generalized Towers of Hanoi problem}
\end{center}
\begin{center}
{\large Mikael Erik Jörgensen}
\end{center}
\begin{center}
{\large Email: mijo0856@student.su.se}
\end{center}

\begin{abstract}
The purpose of this paper is to prove the Frame-Stewart algorithm for the generalized Towers of Hanoi problem as well as finding the number of moves required to solve the problem and studying the multitude of optimal solutions. The main idea is to study how to most effectively move away all but the last disc and use the fact that the total number of moves required to solve the problem is twice this number plus one. 
\end{abstract}
\newpage

{\Large Preface} \newline \newline
It was about four years ago when I first encountered the generalized Towers of Hanoi problem. I had used the original one in teaching purposes and simply extended the problem in a natural way. I also found the Frame-Stewart algorithm but I never proved it since it felt obvious and I assumed the problem was since long proved. Only recently it was pointed out to me that it was in fact not proved, something I could hardly believe. As I searched the web however I found a solution in [Arif] published only about 100 days before. At this point I had already thought of a method of proof which is the one that I present in this paper. Hopefully my ideas may be useful in solving related problems as well.
\newpage

{\large \textbf{Proporties of the puzzle}} \newline
\newline
I shall assume familiarity with the rules of the puzzle. For background information see [wiki] or introduction in [Arif]. \newline
\newline
Now we must label the discs: Disc 1 is the smallest disc, disc 2 the second smallest and so on. \newline
\newline
\textbf{Remark:} A number of words will have specific meaning in this text and are listed in this remark: \newline
$\cdot$ The size of a disc is simply it's assigned number. \newline
$\cdot$ A stack is any number of discs of any sizes located on one and the same peg. \newline
$\cdot$ A tower is a stack such that if disc $i$ is in the tower then either disc $i+1$ is also in the tower or disc $i$ is the largest disc in the tower. \newline
$\cdot$ The height of a stack is the number of discs in it. \newline
$\cdot$ A distribution is a collection of stacks. \newline
\newline
\textbf{Definition 1:} The minimum number of moves required to distribute a single tower consisting of $n$ discs onto $m$ empty pegs with $k$ pegs available is $S(n,k,m)$. We also denote $T(n,k)=S(n,k,1)$. \newline
\newline
Now let's state some useful proporties of the puzzle. \newline
\newline
\textbf{Blocking lemma:} If a puzzle has $k$ pegs and $m$ pegs are occupied by stacks with smallest disc less than $i$ and the i:th disc is on top of it's stack, then the i:th disc may be moved to $k-m-1$ positions. \newline
\newline
{\textbf{Proof}\qquad}
This follows immediately from the rules of the puzzle.
$\textbf{Q.E.D.}$ \newline
\newline 
\textbf{Reversal lemma:} It takes the same number of moves to distribute a stack onto $m$ empty pegs as it takes to reassemble the stack from the distributed position. \newline
\newline
{\textbf{Proof}\qquad}
This is obvious since we can simply move each disc to it's immediate previous position in reverse order and thus it cannot take more moves to reassemble the tower. Now suppose we managed to distribute the tower onto $m$ pegs with a minimal amount of moves i.e. $S(n,k,m)$ moves, then if we could reassemble the tower in fewer moves we could use the reverse strategy to distribute the tower which contradicts the minimality of the distribution strategy.
$\textbf{Q.E.D.}$ \newline
\newline
\textbf{Corollary: Bifurcation theorem:} The number of moves required to solve the general Towers of Hanoi problem is $2S(n-1,k,k-2)+1$. \newline
\newline
{\textbf{Proof}\qquad}
The first thing we need do is to distribute $n-1$ discs onto $k-2$ pegs using $k$ pegs. For this we need exactly $S(n-1,k,k-2)$ moves. Next we move the largest disc (disc n) to it's final position, which requires 1 move. Relabeling the pegs we see that it takes exactly $S(n-1,k,k-2)$ moves to assemble the tower containing all discs onto the destination peg.
$\textbf{Q.E.D.}$ \newline
\newline
\newline
{\large \textbf{Proof of the algorithm}} \newline
\newline
Now we have what we need to prove the Frame-Stewart algorithm in the general case. Notice that all we need to do is to prove that distributing the discs, it is indeed an optimal strategy to construct one tower at a time on the intermediate pegs. \newline
\newline
\textbf{Remark:} When counting the number of moves required to build a stack we count how many moves is made by every disc which will form the stack. \newline
\newline
\textbf{The tower distribution lemma:} $\displaystyle S(n,k,m)=\sum_{i=1}^{m} T(d_i,k-i+1)$ \newline for certain choices $d_i$ where $\displaystyle \sum d_i =n$. \newline
\newline
{\textbf{Proof}\qquad}
Case 1: If $j-1$ stacks are already in place then it is obvious that to form the $j:th$ stack we need $T(d_{j},k-j+1)$ moves. \newline
Case 2: If on the other hand we can construct the $j:th$ stack in less than $T(d_{j},k-j+1)$ then we must first place at least one of the previous stacks on top of the other stacks. Then we can redefine the $d_i$ so that we are back at case 1.
$\textbf{Q.E.D.}$ \newline
\newline
Since the number of moves required to build a stack is independent of the discs sizes we can do this in any way such that it minimizes the number of moves. It is easy to see that in particular we can build only towers resulting in the Frame-Stewart algorithm. \newline
\newline
\textbf{Corollary: The Hanoi theorem:} The Frame-Stewart algorithm is optimal. \newline
\newline \newline \newline
{\large \textbf{Journey numbers} \newline
\newline
For reasons to be seen we will from this point reverse the ordering of the discs so that disc $1$ is the largest disc and disc $2$ is the second largest and so on. [Arif] states that the solution given to the Reve puzzle is in fact unique. This is not true. Even if we think of each peg, except the destination peg, as interchangeable there are many ways to solve the problem and not all are by creating only towers but also other kinds of stacks. An example is the 4 disc Reve puzzle. We may first build a tower of disc 3 and 4. If we do then there are 2 ways to solve the problem depending on wether we move disc 4 or disc 2 first after having moved disc 1. We may also distribute disc 2,3 and 4 onto the 3 available pegs with disc 4 on the destination peg. Then we may choose wether to move disc 4 onto disc 2 or disc 3 and thus there are even 2 different choices for distributing the first 3 discs before moving disc 1. To find the least number of moves required and how many optimal distributions there are we will consider the concept of journey numbers. \newline
\newline
\textbf{Definition 2:} If a stack is built such that each consecutive disc moves the least number of moves possible provided that the stack itself is still built in a minimum number of moves then the disc specific journey number, $j_{k}(i)=$ number of moves made by the $i:th$ largest disc in the stack. \newline
\newline
\textbf{Definition 3:} $\displaystyle J_{k}(n)=\max_{i \le n} j_{k}(i)$ is known as the journey number of a stack with height $n$, or simply the journey number. \newline
\newline
When it is clear from the context the index will be dropped. \newline
\newline
\textbf{Proposition 1:} $j(n) \ge j(n-1)$ \newline
\newline
{\textbf{Proof}\qquad}
Suppose $j(n)<j(n-1)$. Then by the blocking lemma we have $j_{k}(i)=j_{k-1}(i)$ for all $i<n$ which is only possible if $j_{k-1}(i)=2$ for all $1<i<n$. But $j(n)$ cannot be less than $2$ by the blocking lemma which contradicts that $j(n)<j(n-1)$.
$\textbf{Q.E.D.}$ \newline
\newline
\textbf{Corollary 1:} $J(n)=j(n)$ \newline
\newline
{\textbf{Proof}\qquad}
This is obvious since by proposition 1 $\displaystyle \max_{i \le n} j(i)=j(n)$.
$\textbf{Q.E.D.}$ \newline
\newline
\textbf{Corollary 2:} $J_{k}(n)=T(n,k)-T(n-1,k)$ \newline
\newline
{\textbf{Proof}\qquad}
This follows directly from corollary 1.
$\textbf{Q.E.D.}$ \newline
\newline
What we want now is some way to calculate the journey number. If we can find a formula to calculate the journey number we also find a way to calculate $S(n,k,m)$ and as such the number of moves required to solve the general Tower of Hanoi problem. Furthermore we will find how we can distribute the first $n-1$ discs before moving the $n:th$ disc. \newline
\newline
\textbf{Doubling lemma:} Either $J(n)=J(n-1)$ or $J(n)=2J(n-1)$. \newline
\newline
{\textbf{Proof}\qquad}
It is clear that sometimes $J(n)=J(n-1)$ and sometimes not. We need to show that when equality doesn't hold $J(n)=2J(n-1)$. It is clear that $J(n) \le 2J(n-1)$ since we can always stack disc $n$ on to disc $n-1$. Suppose $J(n-1) < J(n) < 2J(n-1)$. Then we must have $J_{k}(i)=J_{k-1}(i)$ for all $i \le n-2$. This is only possible if $J_{k}(i)=2$ for all $2 \le i \le n-2$. But from the bifurcation theorem and corollary 2 we see that $J(n)$ is even and it is easily verified that if $J_{k-1}(n-2)=2$ then $J_{k}(n-1)=2$ which means that $J(n)=3$. A contradiction and we are done.
$\textbf{Q.E.D.}$ \newline
\newline
\textbf{The Journey theorem:} $J(n)=2^r$, whenever \newline $\displaystyle \binom{k+r-3}{k-2} < n \le \binom{k+r-2}{k-2}$. \newline
\newline
{\textbf{Proof}\qquad}
From the doubling lemma we know that $J(n)=2^r$ for some $r$ and as such we only need to prove that $r$ is as stated. It is clear that the statement is true for $J(1)=1$ and $J_{k}(n)=2$ if $k>n>1$. Suppose the statement is true for $J_{k}(m)=2^s$ with $s<r$. By the bifurcation theorem we have $\displaystyle J(n)=2 \max_{i\in I} J_{k-i}(n_i)$ with $2^{r-2}\le J_{k-i}(n_i) \le 2^{r-1}$ where $I=\{0,1,2,\cdots,k-2\}$. Now $\displaystyle n=1+\sum_i n_i$ so we can calculate the maximum possible value of $n$.
$$n=1+\binom{k+r-3}{k-2}+\binom{k+r-4}{k-3}+\cdots+\binom{r-1}{1}=\binom{k+r-2}{k-2}$$
The above calculation is actually easiest seen in pascals triangle, however it is not hard to verify it algebraically either. We have thus established the upper bound for $n$. As for the lower bound we already know it since it must be the upper bound for $r-1$. Thus the induction is complete and the theorem is proven.
$\textbf{Q.E.D.}$ \newline
\newline \newline
{\large \textbf{Consequences of the Journey theorem} \newline
\newline
First consequence of the Journey theorem is that it allows us to calculate the minimum number of moves quite readily. Consider these examples: \newline
\newline
\textbf{Example 1:} Calculate the minimum number of moves for $n=10$ and $k=4$. \newline
\newline
We calculate the journey numbers: $J(1)=1$; $J(s)=2, 2\le s\le3$; $J(s)=4, 4\le s\le6$; $J(s)=8, 7\le s\le10$. We have that it takes $\displaystyle \sum_{i=1}^{10} J(i)=1\cdot 1+2\cdot 2+4\cdot 3+8\cdot 4=1+4+12+32=49$ moves. \newline
\newline
\textbf{Example 2:} Calculate the minimum number of moves for $n=1000$ and $k=30$. \newline
\newline
Same as in Example 1 we calculate journey numbers. $J(1)=1$; $J(s)=2, 2\le s\le29$; $J(s)=4, 30\le s\le15\cdot 29=435$; $J(s)=8, 436\le s\le 4495$. We have that the minimum number of moves is $8\cdot (1000-435)+4\cdot (435-29)+2\cdot (29-2)+1=6317$. Notice that this number is quite a bit smaller that the number of moves it takes to solve the original Towers of Hanoi problem with 13 discs. \newline
\newline
Calculating the number of optimal solutions is hard. Instead of doing this I decided to find a way of calculating the number of main distributions. This is the first step to calculate the actual number of optimal solutions and a necessary one, as far as I can tell. A main distribution is the distribution when all but the largest disc has been moved from the source and the destination peg is empty, i.e. there are no discs on that peg, and we solve the puzzle in the minimum number of moves. We identify any two distributions that has the same stack heights. We shall use journey numbers to find the number of main distributions. \newline
\newline
\textbf{Example 3:} Consider the situation when we have 8 discs and 5 pegs. The journey numbers are $J(1)=1$, $J(2-4)=2$, $J(5-10)=4$. Now if we would have had 10 discs there would have been only one distribution since all the stacks in the main distribution would have had the maximum number of discs for their journey numbers. Now however we have some freedom of placement. We know that in the main distribution we have no stack with journey number greater than 2 and as such we can calculate each stacks maximum height. The first stack cannot be higher than 4 since $J_{5}(5)=4$. The second stack cannot be higher than 3 since $J_{4}(4)=4$ and the last stack cannot be higher than 2 for a similar reason. We can therefor construct stacks of height 4, 2 and 1 or 3, 3 and 1 or 3, 2 and 2. We therefor have 3 main distributions in this case. \newline
\newline
As is clearly noted, journey numbers are crucial in determining the number of main distributions. While this example was easy, the same method can be used in more complicated situations as well. \newline
\newline
\newline
{\large \textbf{Concluding remarks}} \newline
\newline
It is my belief that the tower distribution lemma can be extended to more general Hanoi problems, i.e. problems with the same rules but with different connection between the pegs. In general one may consider a peg as a vertex in a graph. It is my hope that journey numbers may be useful in other Hanoi problems as well. I believe them to be key in finding the number of optimal solutions and the number of moves required to solve such problems. If any of the ideas presented in this paper should come to good use it would make me sincerely happy.

\newpage

{\large REFERENCES} \newline \newline \newline
[Arif] Bijoy Rahman Arif,\newline  On the Footsteps to Generalized Tower of Hanoi Strategy \newline
http://arxiv.org/pdf/1112.0631v1.pdf \newline \newline
[Wiki] Wikipedia article \newline
http://en.wikipedia.org/wiki/Tower\_of\_Hanoi

\end{document}